\documentclass{amsbook}

\theoremstyle{definition}

\theoremstyle{remark}

\def \rr {\mathbb{R}}

\def \rnm {\mathbb{R}^n_-}

\def \huno {H^1_0(\Omega)}

\def \beq {\begin{eqnarray*}}
\def \eeq {\end{eqnarray*}}
\def \beqn {\begin{eqnarray}}
\def \eeqn {\end{eqnarray}}
\def \bequa {\begin{equation}}
\def \eequa {\end{equation}}

\numberwithin{section}{chapter}
\numberwithin{equation}{chapter}

\def\phi{{\varphi}}

\def\cal{{\rm }}

\DeclareSymbolFont{AMSb}{U}{msb}{m}{n}
\DeclareMathSymbol{\N}{\mathbin}{AMSb}{"4E}
\DeclareMathSymbol{\Z}{\mathbin}{AMSb}{"5A}
\DeclareMathSymbol{\R}{\mathbin}{AMSb}{"52}
\DeclareMathSymbol{\Q}{\mathbin}{AMSb}{"51}
\DeclareMathSymbol{\I}{\mathbin}{AMSb}{"49}
\DeclareMathSymbol{\C}{\mathbin}{AMSb}{"43}

\newcommand{ \pOm}{\partial \Omega}

\DeclareSymbolFont{AMSb}{U}{msb}{m}{n}
\DeclareMathSymbol{\N}{\mathbin}{AMSb}{"4E}
\DeclareMathSymbol{\Z}{\mathbin}{AMSb}{"5A}
\DeclareMathSymbol{\R}{\mathbin}{AMSb}{"52}
\DeclareMathSymbol{\Q}{\mathbin}{AMSb}{"51}
\DeclareMathSymbol{\I}{\mathbin}{AMSb}{"49}
\DeclareMathSymbol{\C}{\mathbin}{AMSb}{"43}

\renewcommand{\l }{\lambda }

\renewcommand{\O }{\Omega }
\newcommand{\bequ}{\begin{equation}}
\newcommand{\ee}{\end{equation}}

\def\C{{\mathcal C}}
\def\S{{\mathbb S}}

\def\div{{\rm div}}

\begin{document}
\frontmatter
\title{Functional  Inequalities:\\New Perspectives and New Applications\\ \vskip 20pt
}

\author{Nassif Ghoussoub\footnote{Department of Mathematics, University of British Columbia, Vancouver, B.C. Canada V6T 1Z2. E-mail: nassif@math.ubc.ca. Research partially supported by the Natural Science and Engineering Research Council of Canada.}\qquad
Amir Moradifam\footnote{Department of Mathematics, University of Toronto, Toronto, Canada. E-mail: amir@math.toronto.edu. Research supported by a MITACS postdoctoral fellowship.}\vskip 50pt

\today

}

\maketitle

\setcounter{page}{4}
\mainmatter

%

\chapter*{Preface}

This book is not meant to be another compendium of select inequalities, nor does it claim to contain the latest or the slickest ways of proving them. This project is rather an attempt at  describing how most functional inequalities are not merely the byproduct of ingenious guess work by a few wizards among us, but are often manifestations of certain natural mathematical structures and physical phenomena.  Our main goal here is to show how this point of view leads to ``systematic" approaches for not just proving the most basic functional inequalities, but also for understanding and improving them, and for devising new ones - sometimes at will, and often on demand. 

Our aim is therefore to describe how a few general principles  are behind the validity of  large classes of functional inequalities, old and new. As such, Hardy and Hardy-Rellich type inequalities involving radially symmetric weights are variational manifestations of Sturm's theory on the oscillatory behavior of certain ordinary differential equations.  Similarly, allowable non-radial weights  in Hardy-type inequalities for more general uniformly elliptic operators are closely related to the resolution of certain linear PDEs in divergence form with either a prescribed boundary condition or with prescribed singularities in the interior of the domain.
  
On the other hand, most geometric inequalities including those of Sobolev and Log-Sobolev type, are simply expressions of the convexity of certain free energy functionals along the geodesics of the space of probability measures equipped with the optimal mass transport (Wasserstein) metric.  Hardy-Sobolev and Hardy-Rellich-Sobolev type inequalities are then obtained by interpolating the above inequalities via the classical ones of H\"older. 

Besides leading to new and improved inequalities, these general principles offer new ways for estimating their best constants, and for deciding whether they are attained or not in the appropriate function space.  In Hardy-type inequalities, the best constants are related to the largest parameters for which certain linear ODEs have non-oscillatory solutions. Duality methods, which naturally appear in the new ``geodesic convexity" approach to geometric inequalities,  allow for the evaluation of the best constants from first order equations via the limiting case of Legendre-Fenchel duality, as opposed to the standard method of solving second order Euler-Lagrange equations. 

Whether a ``best constant" on specific domains is attained or not, is often dependent on how it  compares to related best constants on limiting domains, such as the whole space or on half-space. These results are based on delicate blow-up analysis, and are reminiscent of the prescribed curvature problems initiated by Yamabe and Nirenberg.  The exceptional case of the Sobolev inequalities in two dimensions initiated by Trudinger and Moser can also be proved via mass transport methods, and some of their recent improvements by Onofri, Aubin and others are both interesting and still challenging. They will be described in the last part of the monograph.

The part dealing with Hardy and Hardy-type inequalities represents a compendium of work mostly done by --and sometimes with-- my (now former) students Amir Moradifam and Craig Cowan, while the ``mass transport" approach to geometric inequalities follows closely my work with my former student X. Kang and  postdoctoral fellow Martial Agueh. This is largely based on the pioneering work of Cedric Villani, Felix Otto, Robert McCann, Wilfrid Gangbo, Dario Cordero-Erausquin, B. Nazareth, C. Houdr\'e and many others.  The chapters  dealing with Hardy-Sobolev type inequalities follow work done with my students Chaogui Yuan, and Xiaosong Kang, a well as my collaborator Frederic Robert. Finally, much of the progress on the --still unresolved-- best constant in Moser-Onofri-Aubin inequalities on the 2-dimensional sphere was done with my friends and collaborators, Joel Feldman, Richard Froese, Changfeng Gui, and Chang-Shou Lin. I owe all these people a great deal of gratitude.

I thank my wife Louise and my children Mireille, Michelle and Joseph for their patience and their support over the years. 

\aufm{Nassif Ghoussoub}


%




\chapter*{Introduction}

This book is an attempt to describe how a few general principles  are behind the validity of  large classes of functional inequalities, old and new. It consists of six parts, which --though interrelated-- are meant to reflect either the mathematical structure or the physical phenomena behind certain collections of inequalities.  

In Part I, we deal with Hardy-type inequalities involving radially symmetric weights and their improvements. The classical Hardy inequality asserts that for a domain $\Omega$  in $\R^{n}$, $n \geq 3$, with $0 \in \Omega$, the following holds:
\begin{equation}\label{cl-hardy}
\hbox{$\int_{\Omega}|\nabla u |^{2}dx \geq ( \frac{n-2}{2})^{2} 
\int_{\Omega}\frac{u^2}{|x|^{2}}dx$ \quad  for $u \in H^{1}_{0}(\Omega)$.}
\end{equation}
The story here is the newly discovered link between various improvements of this inequality confined to bounded domains and Sturm's theory regarding the oscillatory behavior of certain linear ordinary equations, which we review in Chapter 1. 

In Chapter 2, we first identify suitable conditions on a non-negative  $C^1$-function $P$ defined on an interval $(0, R)$ that will allow for the  following  improved Hardy inequality to hold on every domain $\Omega$ contained in a ball of radius $R$: 
\begin{equation}\label{gen-hardy.0}
\hbox{$\int_{\Omega}|\nabla u |^{2}dx - ( \frac{n-2}{2})^{2} 
\int_{\Omega}\frac{u^{2}}{|x|^{2}}dx\geq \int_{\Omega} P(|x|)u^{2}dx$ \quad for $u \in 
H^{1}_{0}(\Omega)$.}
\end{equation}
It turned out that a necessary and sufficient condition for $P$ to be a {\it Hardy Improving Potential} (abbreviated as {\it HI-potential}) on a ball $B_R$,  is for the following ordinary differential equation associated to $P$ 
\begin{equation}\label{}
y''+\frac{1}{r}y'+P(r)y=0,
\end{equation}
to have a positive solution on the interval $(0, R)$. Elementary examples of HI-potentials are $P \equiv 0$ on any interval $(0,  R)$, $P\equiv 1$ on $(0, z_0)$, where $z_{0}=2.4048...$ is the first zero of the Bessel function $J_0$, and  more generally  $P (r)=r^{-a}$ with $0\leq a<2$ on $(0, z_a)$, where  $z_{a}$ is the first root of the largest solution of the equation $y''+\frac{1}{r}y'+r^{-a}y=0$. Other examples are $P_\rho(r)=\frac{1}{4r^{2}(log\frac{\rho}{r})^{2}}$ on $(0, \frac{\rho}{e})$, but also $P_{k, \rho}(r)=\frac{1}{r^{2}}\sum\limits_{j=1}^k\big(\prod^{j}_{i=1}log^{(i)}\frac{\rho}{r}\big)^{-2}$ on $(0, \frac{\rho}{e^{e^{e^{.^{.^{e(k-times)}}}}}} )$.

Besides leading to a large supply of explicit Hardy improving potentials, this connection to the oscillatory theory of ODEs, gives a new way of characterizing and computing best possible constants such as 
\begin{equation}
\beta (P, R):=\inf_
{\genfrac{}{}{0pt}{}{\scriptstyle{u\in
H^1_0(\O)}}{\scriptstyle{u\neq 0}}}
~\frac{\displaystyle\int_{\O}|\nabla u|^2~dx- \displaystyle \frac{(n-2)^2}{4} \int_{\O}|x|^{-2}|u|^2~dx}{\int_{\O}P(|x|)u^2~dx}.
\end{equation}
On the other hand, the value of the following best constant 
\begin{equation}\mu_{\l}(P, \O):= \inf_
{\genfrac{}{}{0pt}{}{\scriptstyle{u\in
H^1_0(\O)}}{\scriptstyle{u\neq 0}}}
~\frac{\displaystyle\int_{\O}|\nabla u|^2~dx-\l\int_{\O}P(|x|)u^2~dx}
{\displaystyle \int_{\O}|x|^{-2}|u|^2~dx}
\end{equation}
and whether it is attained, depend closely on the position of the singularity  point $0$ vis-a-vis $\Omega$. It is actually equal to  $\frac{(n-2)^2}{4}$, and is never attained in $H^1_0(\Omega)$, whenever $\Omega$ contains $0$ in its interior, but the story is quite different for domains $\Omega$ having $0$ on their boundary. In this case, $\mu_{\l}(P, \O)$ is attained in $H^1_0(\Omega)$ whenever $\mu_{\l}(P, \O)<\frac{n^2}{4}$, which may hold or not. For example, $\mu_{\l}(P, \O)$ is equal to $\frac{n^2}{4}$ for domains that lie on one side of a half-space.

In Chapter 3, we consider conditions on a couple of positive functions  $V$ and $W$ on $(0, \infty)$, which ensure that on some ball $B_R$ of radius $R$ in $\R^{n}$, $n \geq 1$, the 
following inequality holds:
\begin{equation} \label{}
\hbox{$\int_{B}V(|x|)|\nabla u |^{2}dx \geq \int_{B} W(|x|)u^2dx$  \quad for $u \in C_{0}^{\infty}(B_R)$. }
\end{equation}
A necessary and sufficient condition is that the couple $(V, W)$  forms a {\it $n$-dimensional Bessel pair} on the  interval $(0, R)$, meaning that the
equation
\begin{equation}
y''(r)+(\frac{n-1}{r}+\frac{V_r(r)}{V(r)})y'(r)+\frac{W(r)}{V(r)}y(r)=0,
\end{equation}
has a positive solution on $(0, R)$.  This characterization allows us to   improve, extend, and 
unify many results about weighted Hardy-type   inequalities
and their corresponding  best constants. The connection with Chapter 2 stems from the fact that $P$ is a HI-potential  if and only if the couple $(1, \frac{(n-2)^2}{4}r^{-2} +P)$ is a Bessel pair. More generally, the pair 
\begin{equation}
\left(r^{-\lambda}, \, (\frac{n-\lambda-2}{2})^2r^{-\lambda -2}+r^{-\lambda}P(r)\right)
\end{equation}
 is also a $n$-dimensional Bessel pair on $(0, R)$ provided $0\leq \lambda \leq n-2$. Again, the link to Sturm theory provides many more examples of Bessel pairs. 

Hardy's inequality and its various improvements have been used in many contexts such as in the study of the stability of solutions of semi-linear elliptic and parabolic equations, of the asymptotic behavior of the heat equation with singular potentials, as well as in the stability of eigenvalues for 
Schr\"odinger operators. 
In Chapter 4, we focus on applications 
to second order nonlinear elliptic
eigenvalue problems such as 
\begin{equation}
\left\{
\begin{array}{lll}
-\Delta u &=\lambda f(u) &\hbox{in }\Omega \\
\hfill u&=0 &\hbox{on } \pOm,
\end{array}
\right.
\end{equation}
where $ \lambda \ge 0$ is a parameter, $ \Omega$ is a bounded domain in $\rr^N$, $N\geq 2$, and $f$ is a superlinear convex nonlinearity. The bifurcation diagram generally depends on the regularity of the extremal solution, i.e., the one corresponding to the largest parameter for which the equation is solvable. Whether, for a given nonlinearity $f$, this solution is regular or singular depends on the  dimension, and  Hardy-type inequalities are crucial for the identification of the critical dimension. 

Part II deals with the Hardy-Rellich inequalities, which are the fourth order counterpart of Hardy's. In Chapter 5, we show that the same condition on the couple $(V, W)$ (i.e, being a $n$-dimensional Bessel pair)  is also key to improved Hardy-Rellich inequalities of the following type: 
For any radial function $u \in C_{0}^{\infty}(B_R)$ where $B_R$ is a ball of radius $R$  in $\R^{n}$, $n \geq 1$, we have
 \begin{equation} \label{}
\hbox{$\int_{B}V(|x|)|\Delta u |^{2}dx \geq  \int_{B} W(|x|)|\nabla  
u|^{2}dx+(n-1)\int_{B}(\frac{V(|x|)}{|x|^2}-\frac{V_r(|x|)}{|x|})|\nabla u|^2dx$.}
\end{equation}
Moreover, if
\begin{equation}\label{extraordinaire}
\hbox{$W(r)-\frac{2V(r)}{r^2}+\frac{2V_r(r)}{r}-V_{rr}(r)\geq 0$\quad on $[0, R),$} 
\end{equation}
then the above inequality holds true for all $u \in C_{0}^{\infty}(B_R)$ and not just the radial ones.  By combining this with the inequalities involving the Dirichlet integrals of Chapter 3, one obtains various improvements of the Hardy-Rellich inequality for $H^2_0(\Omega)$. In particular, for any bounded domain $\Omega$ containing $0$ with $\Omega \subset B_R$, we have the following inequality for all $u \in H^2_0(\Omega)$,
\begin{equation}
\int_{\Omega}|\Delta u|^{2}dx \geq 
\frac{n^{2}(n-4)^2}{16}\int_{\Omega}\frac{u^2}{|x|^4}dx+\frac{\beta (P;  R)\big(n^2+(n-\lambda-2)^2\big)}{4}
\int_{\Omega}\frac{P(|x|)}{|x|^2}u^2 dx,
\end{equation}
where $n\geq 4$, $\lambda <n-2$, and where $P$ is a HI-potential on $(0, R)$ such that $\frac{P_r(r)}{P(r)}=\frac{\lambda}{r}+f(r)$, $f(r)\geq 0$ and
$\lim\limits_{r \rightarrow 0}rf(r)=0$.

In Chapter 6, we explore Hardy-type inequalities for $H^1(\Omega)$-functions, i.e., for functions which do not necessarily have compact support in $\Omega$. In this case, a penalizing term appears in order to account for the boundary contribution. If a pair of positive radial 
functions  $(V, W)$ is a n-dimensional Bessel pair on an interval $(0, R)$, and if $B_R$ is a ball of radius $R$  in $\R^{n}$, $n
\geq 1$, then there exists $\theta >0$ such that the following inequality holds:
\begin{equation} \label{}
\hbox{$\int_{B_R}V(x)|\nabla u |^{2}dx \geq  \int_{B_R} W(x) 
u^{2}dx-\theta\int_{\partial B_R}u^2 ds$\ \ for\, $u \in H^1(B_R)$,}
\end{equation}
and for all radial functions $u \in H^2(B_R)$,
\begin{eqnarray}
\qquad \quad \int_{B_R}V(|x|)|\Delta u |^{2}dx &\geq&  \int_{B_R} W(|x|)|\nabla
u|^{2}dx+(n-1)\int_{B_R}(\frac{V(|x|)}{|x|^2}-\frac{V_r(|x|)}{|x|})|\nabla
u|^2dx \\ && +\big[(n-1)-\theta) V(R)\big]\int_{\partial B_R}|\nabla u|^2\, dx.\nonumber
\end{eqnarray}
The latter inequality holds for all functions in $H^2(B)$ provided (\ref{extraordinaire}) holds.
The combination of the two inequalities lead to various weighted Hardy-Rellich
inequalities on $H^{2}\cap H^{1}_{0}$.

In Chapter 7,  we investigate some applications of the improved
Hardy-Rellich inequalities to fourth order nonlinear elliptic
eigenvalue problems of the form
\begin{equation}
\left\{
\begin{array}{ll}
\Delta^2 u = \lambda f(u) & \hbox{in }\Omega \\
u =\Delta u = 0 &\hbox{on } \pOm,
\end{array}
\right.
\end{equation}
as well as their counterpart with Dirichlet boundary conditions. In particular, they are again crucial for the identification of ``critical dimensions" for such equations involving either an exponential or a singular supercritical nonlinearity.

Part III addresses  Hardy-type inequalities for more general uniformly elliptic operators.  The issue of allowable non-radial weights (to replace $\frac{1}{|x|^2}$) is then closely related to the resolution of certain linear PDEs in divergence form with either prescribed conditions on the boundary or with prescribed singularity in the interior.  We also include $L^p$-analogs of various Hardy-type inequalities.

In Chapter 8, the following general Hardy inequality is associated to any given symmetric, uniformly positive definite $ n \times n$ matrix $ A(x)$ defined in $ \Omega$ with the notation $ | \xi |_A^2:=\langle A(x) \xi ,  \xi\rangle $ for $ \xi \in \R^n$.
\begin{equation}\label{hardy.100}
  \int_\Omega | \nabla u|_A^2 dx \ge \frac{1}{4} \int_\Omega \frac{| \nabla E|_A^2}{E^2}u^2dx, \qquad u \in H_0^1(\Omega)
\end{equation}
   The basic assumption here is that $ E$ is a positive solution to $ -\div(A \nabla E)\, dx=\mu$ on $\Omega$, where $\mu$ is any nonnegative nonzero finite measure on $ \Omega$.  The above inequality 
      is then optimal in either one of  the following  two cases:
      \begin{itemize}
        \item {\it $E$ is an interior weight}, that is $ E=+\infty$ on the support of $ \mu$, or 
    \item   {\it $E$ is a boundary weight}, meaning that  $ E=0$ on $ \pOm$.
           \end{itemize}
    While the case of an interior weight extends the classical Hardy inequality, the case of a boundary weight extends the following so-called {\it Hardy's boundary inequality}, which holds for any bounded convex domain $ \Omega \subset \R^n$ with smooth boundary:
 \begin{equation} \label{} 
\hbox{$ \int_\Omega | \nabla u|^2 dx \ge \frac{1}{4} \int_\Omega \frac{u^2}{{\rm dist}(x,\pOm)^2} dx$ for $ u \in H_0^1(\Omega)$. }
 \end{equation}     
 Moreover the constant $ \frac{1}{4}$ is optimal and not attained.  One also obtains other Hardy inequalities involving more general distance functions.  For example, if $ \Omega$ is a domain in $ \R^n$ and $M$ a piecewise smooth surface of co-dimension $k$ ($ k=1, ... , n$). 
 Setting $ d(x):= {\rm dist}(x,M)$ and suppose $ k \neq 2$ and $ - \Delta d^{2-k} \ge 0 $ in $ \Omega \backslash M$,  then 
 \begin{equation} \label{AAA}
\hbox{$ \int_\Omega | \nabla u|^2 dx \ge \frac{(k-2)^2}{4} \int_\Omega \frac{u(x)^2}{d(x)^2} dx$ \quad for $ u \in H_0^1( \Omega \backslash M)$.  }
 \end{equation}
The inequality is not attained in either case, and  one can therefore get  the following improvement for  the case of a boundary weight:   
\begin{equation} \label{}  \int_\Omega | \nabla u|_A^2dx \ge \frac{1}{4} \int_\Omega \frac{| \nabla E|_A^2}{E^2}u^2dx + \frac{1}{2} \int_\Omega \frac{u^2}{E} d \mu, \qquad u \in H_0^1(\Omega) 
\end{equation} 
which is optimal and still not attained.  
 Optimal weighted versions of these inequalities are also established, as well as  their $L_p$-counterparts  when $ p \neq 2$.  
 Many of the Hardy inequalities obtained in the previous chapters 
 can be recovered via the above approach, by using  suitable choices for $E$ and $A(x)$. 
 
   In Chapter 9,  we investigate the possibility of improving (\ref{hardy.100}) in the spirit of Chapters 3 and 4, namely whether one can find conditions on non-negative potentials $V$ so that the following improved inequality holds:
 \begin{equation} \label{THREE.100}
  \hbox{$  \int_\Omega | \nabla u|_A^2dx - \frac{1}{4} \int_\Omega \frac{| \nabla E|_A^2}{E^2} u^2dx \ge \int_\Omega V(x) u^2dx$ for $ u \in H_0^1( \Omega)$.} 
    \end{equation}  
 Necessary and sufficient conditions on $V$ are given for (\ref{THREE.100}) to hold, in terms of the solvability of a corresponding linear PDE.  Analogous results involving improvements are obtained for the weighted versions.  Optimal inequalities are also obtained for 
 $ H^1(\Omega)$.

 We conclude Part III by considering in Chapter 10, applications of the Hardy inequality for general uniformly elliptic operators to study the regularity of stable solutions of certain nonlinear eigenvalue problems involving advection such as 
\begin{equation} 
\left\{
\begin{array}{lll}
-\Delta u + c(x) \cdot \nabla u &=& \frac{\lambda}{(1-u)^2} \quad \mbox{in $ \Omega$}, \\
\hfill u &=& 0 \qquad  \mbox{on  $ \pOm$}, \\
\end{array}
\right.
\end{equation} where $ c(x)$ is a smooth bounded vector field on $\bar \Omega$. 

  In Part IV, we describe how the Monge-Kantorovich theory of mass transport provides a framework that encompasses most geometric inequalities. Of importance is the concept of {\it relative energy of $\rho_0$ with
respect to
$\rho_1$} defined as:
  \begin{equation} {\rm
H}^{F,W}_{V}(\rho_0|\rho_1):={\rm H}^{F,W}_{V}(\rho_0)- {\rm
H}^{F,W}_{V}(\rho_1),
\end{equation}
where $\rho_0$ and $\rho_1$ are two probability densities, and where the {\it Free Energy Functional} ${\rm H}_V^{F,W}$ is defined on the set 
${\mathcal P}_a(\Omega)$ of probability densities on a domain $\Omega$ as:
\begin{equation}
{\rm H}_V^{F,W}(\rho):=\int_{\Omega}\left[F(\rho) +\rho V
+\frac{1}{2}(W\star\rho)
\rho\right]\,\mbox{d}x.
\end{equation}
${\rm H}_V^{F,W}$ being the sum of the internal energy ${\rm H^F}(\rho):=\int_{\Omega}F(\rho) dx$, 
the potential energy ${\rm H}_V(\rho):=\int_{\Omega} \rho V dx$
and the interaction energy ${\rm H}^{W}(\rho):=\frac{1}{2}\int_{\Omega}\rho(W\star\rho)\,\mbox{d}x$. Here  $F$ is a differentiable function on $(0,\infty)$, while the confinement (resp., interactive) potential $V$ (resp., $W$) are $C^2$-functions on $\R^n$ satisfying $D^2V\geq \mu I$ (resp.,   $D^2W\geq \nu I$) for some $\mu, \nu \in \R$.

In Chapter 11, we describe Brenier's solution of the Monge problem with quadratic cost, which yields that the  
Wasserstein distance $W(\rho_0, \rho_1)$ between two probability densities $\rho_0$ on $X$  and $\rho_1$ on $Y$, i.e., 
\begin{equation}\label{Wass}
W(\rho_0, \rho_1)^2=\inf \left\{\int_{X}{|x-s(x)|^2}dx;\,  s\in {\cal S}(\rho_0, \rho_1) \right\}
\end{equation}
is achieved by  the gradient $\nabla \phi$ of a convex function $\phi$. 
Here $S(\rho_0, \rho_1)$ is the class of all Borel measurable maps $s: X\to Y$  that ``push" $\rho_0$ into $\rho_1$, i.e., those which satisfy the change of variables formula,
\begin{equation}
\int_{Y}{h(y)\rho_1(y)}dy=\int_{X}{h(s(x))\rho_0(x)}dx\ \ \ \text{for every}\ \ h \in C(Y).
\end{equation}
This fundamental result allows one to show that for certain natural candidates $F, V$ and $W$, the corresponding free energy functionals ${\rm H}_V^{F,W}$ are convex on the geodesics of optimal mass transport joining two probability densities in ${\mathcal P}_a(\Omega)$. This  convexity property  translates into a very general inequality relating the relative total energy between the initial and final configurations $\rho_0$ and $\rho_1$, to their entropy production ${\mathcal I}_{c^*}(\rho|\rho_{_V})$, their Wasserstein distance $W_2^2(\rho_0, \rho_1)$, as well as the Euclidean distance between their barycenters $|{\rm b}(\rho_0)-{\rm b}(\rho_1)|$,
\begin{equation}
\label{ultimate}
{\rm H}^{^{F,W}}_{_{V+c}}(\rho_0|\rho_1)+\frac{\lambda+\nu}{2}
W_2^2(\rho_0,
\rho_1)-\frac{\nu}{2}|{\rm b}(\rho_0)-{\rm b}(\rho_1)|^2
  \leq {\rm H}_{c+\nabla V\cdot x}^{^{-nP_F,2x\cdot\nabla W}}(\rho_0) +
{\cal I}_{c^*}(\rho|\rho_{_V}).
\end{equation}
Here  $P_F(x):=xF^\prime(x)-F(x)$ is the {\it pressure function} associated to $F$, while $c$ is a Young function (such as $c(x)=\frac{1}{p}|x|^p$), $c^*$ is its Legendre transform, while ${\cal I}_{c^*}(\rho|\rho_{_V})$ is the 
{\it relative entropy
production-type function of $\rho$ measured
against $c^*$} defined as
\begin{equation}
{\cal I}_{c^*}(\rho|\rho_{_V}):=\int_\Omega \rho
c^\star\left(-\nabla\left(F^\prime(\rho)+V+W\star\rho\right)\right)\,\mbox{d}x.
\end{equation}
Once this general comparison principle is established, various -- new and old -- inequalities follow by simply considering different examples of internal, potential and interactive energies, such as $F(\rho)=\rho \ln \rho$ or $F(\rho)=\rho^\gamma$, and $V$ and $W$ are convex functions (e.g., $V(x)=\frac{1}{2}|x|^2$), while $W$ is required to be even.    

The framework is remarkably encompassing even when $V=W\equiv 0$, as it is shown in Chapter 12 that the following inequality, which relates the internal energy of a probability density $\rho$ on $\R^n$ to the corresponding entropy production contains almost all known Euclidean Sobolev and log-Sobolev inequalities:
  \begin{equation}
\label{eqn:7} \int_{\Omega}[F(\rho)+ n P_F(\rho)]\, dx \leq 
\int_\Omega\rho c^\star\left(-\nabla
(F^\prime\circ\rho)\right)\,\mbox{d}x+K_c.
\end{equation}
The latter constant $K_c$ can always be computed from $F$ and the Young function $c$. 

The approach allows for a direct and unified way for computing best constants and extremals.  It also leads to remarkable duality formulae, such as the following associated to the standard Sobolev inequality for $n\geq 3$ and where $2^*:=\frac{2n}{n-2}$:
  \begin{eqnarray}
\label{}
 \qquad \sup\big\{\frac{n(n-2)}{n-1}\int_{\R^n}\rho (x)^{\frac{n-1}{n}}dx
- \int_{\R^n}|x|^2\rho (x)dx; \int_{\R^n}\rho(x)\ dx=1\big\} \qquad  \qquad   \\
\qquad  \qquad  \qquad \qquad =
\inf\Big\{\int_{\R^n} |\nabla f|^2dx;\
f\in C^\infty_0(\R^n);
\
\int_{\R^n}|f|^{2^*}dx=1\Big\}. \nonumber
\end{eqnarray}
This type of duality also yields a remarkable correspondence between ground state solutions of certain
quasilinear (or semi-linear) equations, such as  ``Yamabe's'', 
\[
-\Delta f =  |f|^{2^*-2}f \ {\rm on }\ \R^{n},
\]
 and stationary solutions of the (non-linear) Fokker-Planck equations 
$                  \frac{\partial u}{\partial t}= \Delta u^{1-{1\over n}} +
{\rm div}(x.u)$,  which --after appropriate scaling-- reduces to the fast diffusion equation 
\[ \frac{\partial u}{\partial t}= \Delta u^{1-{1\over n}} \ {\rm on }\ \R^+\times\R^{n}.
\]

Chapter 13 deals with applications to Gaussian geometric inequalities.  We first establish the so-called HWBI inequality, which follows immediately from a direct application of (\ref{ultimate}) with paramatrized quadratic Young functions $c_\sigma(x)=\frac{1}{2\sigma}{|\,x\,|^2}$ for
$\sigma>0$, coupled with a simple scaling argument:  
\begin{equation}
\label{}
{\rm H}^{F, W}_V(\rho_0|\rho_1)\leq
W_2(\rho_0,\rho_1)\sqrt{I_2(\rho_0|\rho_V)}
-\frac{\mu+\nu}{2}W_2^2(\rho_0,\rho_1)+\frac{\nu}{2}|{\rm b}(\rho_0)-{\rm
b}(\rho_1)|^2.
\end{equation}
This gives a unified approach for --extensions of-- various
powerful
inequalities by Gross, 
Bakry-Emery,  
Talagrand, 
Otto-Villani, 
Cordero-Erausquin,
and others. As expected, such
inequalities also lead to exponential rates of convergence to equilibria
for solutions
of Fokker-Planck and McKean-Vlasov type equations. 

Part V deals with Caffarelli-Kohn-Nirenberg and Hardy-Rellich-Sobolev type inequalities. All these can be  obtained by simply interpolating --via H\"older's inequalities-- many of the 
previously obtained inequalities. This  is done in Chapter 14, where it is also shown that the best constant  in the Hardy-Sobolev inequality, i.e., 
\begin{equation}
\mu_s (\Omega):= \inf \left \{ \int_{\Omega}| \nabla u|^2
dx;\, u \in \huno \hbox{ and }  \int_{\Omega} \frac {|u|^{2^*(s)}}{|x|^s}\,
dx =1\right\},
\end{equation}
where $0<s<2$ and $2^*(s)=\frac{2(n-s)}{n-2}$,  is never attained when $0$ is in the interior of the domain $\Omega$, unless the latter is the whole space $\R^n$, in which case explicit extremals are given. This is not the case when $\Omega$ is half-space $\rnm$, where only the symmetry of the extremals is  shown. Much less is known about the extremals in the Hardy-Rellich-Sobolev inequality (i.e., when $s>0$) even when $\Omega=\R^n$. 

The problem whether $\mu_s(\Omega)$ is attained becomes more interesting when $0$ is on the 
boundary $\partial \Omega$ of the domain $\Omega$.  The attainability is then closely related to the geometry of
$\partial\Omega$, as we show in chapter 15, that the 
negativity of the mean curvature   of $\partial \Omega $ at $0$ is sufficient to ensure the
attainability of $\mu_{s}(\Omega)$.

In Chapter 16, we consider log-Sobolev inequalities on the line, such as those involving the functional  
\begin{equation}
I_\alpha (g) = {\alpha\over 2}\int_{-1}^1 (1-x^2)|g'(x)|^2\ dx +
\int_{-1}^1 g(x)\ dx - \ln{ 1\over 2 }\int_{-1}^1 e^{2g(x)}\ dx
\end{equation}
on the space $H^1(-1,1)$ of $L^2$-functions on $(-1,1)$ such that
 $
(\int_{-1}^1(1-x^2)|g'(x)|^2 dx)^{1/2} <\infty
$.  
We then show that  if $J_\alpha$ is restricted to the manifold
\[
{\mathcal G}=\left\{g \in H^{1}(-1,1); \,   \int_{-1}^1 e^{2g(x)} x dx = 0\right\},
\]
then the following hold:
\begin{equation}\label{ffgg.99}
\hbox{$\inf\limits_{g\in {\mathcal G}} I_\alpha(g) = 0$\,\,   if $\alpha \ge {1 \over 2}$, \quad and \quad $\inf\limits_{g\in {\mathcal G}} I_\alpha(g) = -\infty$\,\,  if $\alpha < {1\over 2}$.}
 \end{equation}
We also give a recent result of Ghigi, which says that the functional \begin{equation}
  \Phi(u) = \int_{-1}^1 u(x) \, d  x - \log \biggl (  \frac {1} {2} \int_{-\infty}^{+\infty} e^{-2u^*(x)}\, d  x \biggr )
\end{equation}
is convex  on the cone $\mathcal W$ of all bounded convex functions $u$ on $(-1,1)$, where here $u^*$ denotes the Legendre transform of $u$, and that 
\[\inf_{u\in \mathcal W} \Phi (u)=\log (\frac{4}{\pi}).\]
Both inequalities play a key role in the next two chapters, which address  inequalities on the two-dimensional sphere $\S^2$. It is worth noting that Ghigi's inequality relies on the Pr\'ekopa-Leindler principle, which itself is another manifestation of a mass transport context. One can therefore infer that the approach of Part IV can and should be made more readily applicable to critical Moser-type inequalities. 

In Chapter 17, we establish the Moser-Trudinger inequality, which states that for $\alpha \geq 1$, the functional 
\begin{equation}
J_\alpha(u): =  \alpha\int_{S^2}|\nabla u|^2\, d\omega
+ 2 \int_{\S^2} u\, d\omega
-\ln \int_{\S^2} e^{2u}\, d\omega
\end{equation}
is bounded below on the  Sobolev space $H^{1,2}(\S^2)$, where here $d\omega:=\frac{1}{4\pi}\sin
\theta\, d\theta \wedge d\phi$ denotes Lebesgue measure on the unit sphere, normalized so that $\int_{\S^2} d\omega = 1$.  We also give a proof of Onofri's inequality which states that the infimum of $J_\alpha$ on $H^{1,2}(\S^2)$ is actually equal to zero for all $\alpha \geq 1$, and that 
\begin{equation}
\inf\{J_1 (u);\, u\in H^{1,2}(\S^2) \}=\inf\{J_1 (u);\,  u\in {\mathcal M}\}= 0,
\end{equation}
where ${\mathcal M}$ is the submanifold 
${\mathcal M}=\{u\in H^{1,2}(\S^2);\,  \int_{\S^2}e^{2u}{\bf x} d\omega = 0\}.$ Note that  this inequality, once  applied to axially symmetric functions, leads to the following counterpart of (\ref{ffgg.99})
\begin{equation}
\hbox{$\inf\limits_{g\in H^1(-1,1)}I_\alpha(g)=\inf\limits_{g\in {\mathcal G}} I_\alpha(g)=0$\qquad  if $\alpha \geq 1$.}
\end{equation}
In Chapter 18, we address results of T. Aubin asserting that once restricted to the submanifold ${\mathcal M}$, the functional $J_\alpha$ then remains bounded below (and coercive) for smaller values of $\alpha$, which was later conjectured by A. Chang and P. Yang to be equal to $\frac{1}{2}$. We conclude the latest developments on this conjecture,  including a proof that 
\begin{equation}
\hbox{$\inf\{J_{\alpha} (u);\, u\in {\mathcal M}\}= 0$ \quad if $\alpha \geq  \frac{2}{3}$.}
\end{equation}
The conjecture remains open for $1/2 <\alpha <2/3$.

We have tried to make this monograph as self-contained as possible. That was not possible though, when dealing with the applications such as in Chapters 4, 7 and 10. We do however give enough references for the missing proofs. 

The rapid development of this area and the variety of applications forced us to be quite selective. We mostly concentrate on certain recent advances not covered in the classical books such as the one by R. A. Adams \cite{Ad} and V. G. Maz'ya \cite{M}.  Our choices reflect our taste and what we know --of course-- but also our perceptions of what are the most fundamental functional inequalities, the novel methods and ideas, those that are minimally ad-hoc, as well as the ones we found useful in our work. 
It is however evident that this compendium is far from being an exhaustive account of this continuously and rapidly evolving line of research. One example that comes to mind are inequalities obtained by  interpolating between the Hardy and the Trudinger-Moser inequalities. One then gets the singular Moser-type inequalities, which states that for some $C_{0}=C_{0}(n,\left\vert
\Omega\right\vert )>0$,  one has for any $u\in W_{0}^{1,n}\left(
\Omega\right)  $ with $\int_{\Omega}\left\vert \nabla u\right\vert
^{n}dx\leq1$,
\begin{equation}
\int_{\Omega}\frac{\exp\left(  \beta\left\vert u\right\vert ^{\frac{n}{n-1}%
}\right)  }{\left\vert x\right\vert ^{\alpha}}dx\leq C_{0}, 
\end{equation}
for any $\alpha\in\left[  0,n\right)  ,~0\leq\beta\leq\left(
1-\frac{\alpha}{n}\right)n\omega_{n-1}^{\frac{1}{n-1}}$, where
$\omega_{n-1}=\frac{2\pi^{\frac{n}{2}}}{\Gamma\left(  \frac{n}{2}\right)  }$
is the area of the surface of the unit $n$-dimensional ball. See for instance \cite{Ad2, AdY}).

 The recent developments on these inequalities could have easily constituted a Part VII for this book, but we had to stop somewhere and this is where we stop.



\backmatter

%

\bibliographystyle{amsalpha}

\end{document}